\def\arctanh{{\rm arctanh}}

\def\zbar{\overline{z}}
\def\tr{{\rm tr}}

\def\IH{{\mathbb H}^3}

\def\IR{{\mathbb R}}

\def\II{{\mathbb I}}

\def\ID{{\mathbb D}}

\def\oR{\overline {\IR}}
\def\IS{{\mathbb S}}
\def\IZ{{\mathbb Z}}
\def\IH{{\mathbb H}}
\def\IC{\mathbb C}
\def\IQ{\mathbb Q}
\def\oC{\hat{\IC}}

\def\zbar{{\overline{z}}}
\def\hbar{{\overline{h}}}
\def\zetabar{{\overline{\zeta}}}

\documentclass{article}

\usepackage{amsthm, amssymb}
\usepackage{amsmath}
\usepackage{amsfonts}
\usepackage{epsfig,multicol}
\usepackage{pstricks}
\usepackage{graphicx}

\title{{\bf Random Kleinian Groups, {\bf II} }\\ {Two parabolic generators}}

\author{Gaven Martin, Graeme O'Brien and Yasushi Yamashita \thanks{Research supported in
part by grants from the N.Z.
Marsden Fund.  \newline \newline AMS
(1991) Classification.
Primary 30C60, 30F40, 30D50, 20H10, 22E40, 53A35, 57N13, 57M60}  }

\date{}

\begin{document}

\maketitle
\newtheorem{theorem}{Theorem}[section]    
                                           
\newtheorem{lemma}[theorem]{Lemma}         
\newtheorem{corollary}[theorem]{Corollary} 
\newtheorem{remark}[theorem]{Remark}       
\newtheorem{definition}[theorem]{Definition}
\newtheorem{conjecture}[theorem]{Conjecture}
\newtheorem{proposition}[theorem]{Proposition}
\newtheorem{example}[theorem]{Example}

\newcommand{\param}{(\gamma,\beta,\beta')}
\newcommand{\parfour}{(\gamma,\beta,-4)}

\numberwithin{equation}{section}

\newcommand{\abs}[1]{\lvert#1\rvert}

\renewcommand{\theequation}{\thetheorem} 
                    
\makeatletter
\let \c@equation=\c@theorem

\begin{abstract}  In earlier work we introduced geometrically natural probability measures on the group of all M\"obius transformations in order to study ``random'' groups of M\"obius transformations,  random surfaces,  and in particular random two-generator groups,  that is groups where the generators are selected randomly, with a view to estimating the likely-hood that such groups are discrete and then to make calculations of the expectation of their associated parameters,  geometry and topology.  In this paper we continue that study and identify the precise probability that a Fuchsian group generated by two parabolic M\"obius transformations is discrete,  and give estimates for the case of Kleinian groups generated  by a pair of random parabolic elements which we support with a computational investigation into of the Riley slice as identified by Bowditch's condition,  and establish rigorous bounds.
\end{abstract}

\maketitle

\section{Introduction.}   There is a consider literature on the question of discreteness of Kleinian groups generated by two parabolic elements.  This starts with the work of Shimitzu and Leutbecher \cite{Le},  and followed by Lyndon and Ullman \cite{LU},  through to the seminal work of Riley \cite{Riley} who clarified the important connections with hyperbolic two bridge knot and link complements after the work of Thurston.  There is much in between.   Stunning visualisations of aspects of this work can be found in \cite{Indra}.

\medskip

In this paper we recall the notion of a random Kleinian group and study the case of two generator groups conditioned by the assumption that both generators are parabolic. In this setting we are able to give the precise probability that such a group is discrete in the Fuchsian case,  and reasonable estimates in the case of Kleinian groups.  We should expect that almost surely (that is with probability one) a finitely generated subgroup of the M\"obius group is free.   We then discuss alternate notions of random groups in terms of probability measures on moduli spaces and give a computationally supported calculation of the area of the Riley slice based on Bowditch's condition \cite{Bow} so as to find another estimate for the probability of a random  group being discrete,  given it has parabolic generators.

\section{Random Fuchsian Groups.}  

We first give definitions in the context of Fuchsian groups to support our later definitions and results. 

\medskip

If $A\in PSL(2,\IC)$ has the form
\begin{equation}\label{Fspace}
A =\pm \left( \begin{array}{cc} a & c \\ \bar c & \bar a \end{array} \right), \hskip15pt |a|^2-|c|^2 = 1,
\end{equation}
then the associated linear fractional transformation $f:\oC\to\oC$ defined by
\begin{equation}\label{fdef}
f(z) = \frac{az + c}{\bar c z + \bar a}
\end{equation}
preserves the unit circle as $\left| \frac{az + c}{\bar c z + \bar a} \right| = |\zbar| \left| \frac{az + c}{\bar a \zbar+\bar c |z|^2} \right|$.
The rotation subgroup ${\bf K}$ of the disk,  $z\mapsto \zeta^2 z$, $|\zeta|=1$,   and the nilpotent or parabolic subgroup (conjugate to the translations) have the respective representations
\[  \left( \begin{array}{cc} \zeta & 0 \\ 0 & \bar \zeta \end{array} \right), \;\;\;|\zeta|=1, \hskip15pt   \left( \begin{array}{cc} 1+it& t \\ t & 1-i t\end{array} \right),\;\;\; t\in \IR . \]
The group of all matrices satisfying (\ref{Fspace}) will be denoted ${\cal F}$ and we refer to ${\cal F}$ as Fuchsian space.  It is not difficult to construct an algebraic isomorphism ${\cal F}\equiv PSL(2,\IR)\equiv Isom^+(\IH^2)$,  the isometry group of two-dimensional hyperbolic space, and we will often abuse notation by moving between $A$ and $f$ interchangeably.   We also seek probability distributions from which we can make explicit calculations and are geometrically natural (see in particular Lemma \ref{3.2}) and therefore impose the following distributions on the entries of this space of matrices ${\cal F}$.  We select 
\begin{itemize}
\item(i) $\zeta=a/|a|$ and $\eta=c/|c|$  uniformly in the circle $\IS$, with arclength measure,    and 
\item(ii) $t=|a|\geq 1$ so that
\[ 2\arcsin(1/t) \in [0,\pi] \]
is uniformly distributed.    
\end{itemize} 
Notice that both the product  $\zeta\eta$ is uniformly distributed on the circle as a simple consequence of the rotational invariance of arclength measure.  Further,  this measure is equivalent to the uniform probability measure $\arg(a)\in [0,2\pi]$. Next observe the pdf for $|a|$.
\begin{lemma}  The random variable $|a|\in [1,\infty)$ has the pdf
\[ F_{|a|}(x)= \frac{2}{\pi} \; \frac{1}{x\sqrt{x^2-1}} \]
\end{lemma}
Notice that the equation $1+|c|^2=|a|^2$ tells us that $\arctan(\frac{1}{|c|})$ is uniformly distributed in $[0,\pi]$.  Every M\"obius transformation of the unit disk $\ID$ can be written in the form
 \begin{equation}\label{mob}
 z\mapsto \zeta^2 \,  \frac{z-w}{1-\bar w z}, \hskip10pt |\zeta|=1, w\in \ID
 \end{equation}
The matrix representation of (\ref{mob}) in the form (\ref{Fspace}) is 
 \[ \zeta^2 \,  \frac{z-w}{1-\bar w z} \leftrightarrow \left(\begin{array}{cc} \frac{\zeta}{\sqrt{1-|w|^2} } &  -\frac{\zeta w}{\sqrt{1-|w|^2}} \\- \frac{\zetabar \bar w}{\sqrt{1-|w|^2} }& \frac{\zetabar}{\sqrt{1-|w|^2}  }
   \end{array} \right), \]
so $\zeta$ and $\frac{w}{|w|}$ are uniformly distributed in $\IS$ and $ \arccos(|w|)=\arcsin(\sqrt{1-|w|^2} ) \in [0,\pi/2]$
is uniformly distributed,  with   $|w|$ having the p.d.f. $\frac{4}{\pi} \sqrt{1-y^2}$,  $y\in [0,1]$. 

Next (\ref{mob}) represents a parabolic transformation when $|w|=|\sin(\theta)|$,  $\theta=\arg(\zeta)$.  Thus we observe
\begin{lemma}  Let $f\in {\cal F}$ be a random parabolic element.  Then the expected value of $|f(0)|$ is $\frac{4}{3\pi}\approx 0.4244\ldots$ with variance $\sigma^2=\frac{1}{4}-\frac{16}{9 \pi ^2}\approx 0.0698\ldots$.
\end{lemma}

\medskip

The isometric circles of the M\"obius transformation $f$ defined at (\ref{fdef}) are the two circles
\[  C_+ =\Big\{|z+\frac{\bar a}{\bar c} | = \frac{1}{|c|}\Big\},  \hskip10pt C_-=\Big\{z:|z-\frac{a}{\bar c}|=\frac{1}{|c|} \Big\} \]
which are paired by the action of $f$ and $f^{-1}$,   $f^{\pm1}(C_{\pm})=C_{\mp}$.  The {\em isometric disks} are the finite regions bounded by these two circles.  Since $|a|^2=1+|c|^2\geq 1 $,  both these circles meet the unit circle in an arc of angle $\theta\in [0,\pi]$.  Some elementary trigonometry reveals that
$\sin \frac{\theta}{2} = \frac{1}{|a|}
$. Thus by our choice of distribution for $|a|$ we obtain the following key result to justify our assertion of being geometrically natural.

\begin{lemma} \label{3.2} The arcs determined by the intersections of the finite disks bounded by the the isometric circles of $f$,  where $f$ is chosen according to the distribution (i) and (ii),  are centred on uniformly distributed points of $\IS$ and have arc length uniformly distributed in $[0,\pi]$.  
\end{lemma}

\section{Fuchsian groups generated by two parabolics.}  If  $f$ and $g$ are parabolic elements of ${\cal F}$ then we have the matrix representations
\begin{equation}\label{2para} f \leftrightarrow \pm \left( \begin{array}{cc} 1+ix & x e^{i\theta} \\ x e^{-i\theta} & 1-i x \end{array} \right), g \leftrightarrow \pm \left( \begin{array}{cc} 1+iy & y e^{i\psi} \\y e^{-i\psi} & 1-i y \end{array} \right),\hskip15pt  \theta,\psi \in_u [0,2\pi] \end{equation}
and $x=\cot(\eta)$ with $\eta$ chosen uniformly from $[0,\pi]$,   $\eta\in_u[0,\pi]$.  The important invariant - and that which determines this group uniquely up to conjugacy in $PSL(2,\IC)$ is the commutator parameter
\begin{equation}
\gamma = \gamma(f,g) = \tr [f,g]-2, \hskip20pt [f,g]=fgf^{-1}g^{-1}.
\end{equation}
Given the above matrix representations at (\ref{2para}) we obtain
\begin{equation}
\gamma(f,g) = 16 x^2 y^2 \sin^4 \big(\frac{\theta-\psi}{2}\big)
\end{equation}
Notice that $\theta-\psi = \arg(e^{i\theta}e^{i\psi})$ and rotational invariance of arclength measure shows $\theta-\psi$ to be uniformly distributed and so therefore is $\eta=(\theta-\psi)/2$.  Thus to consider the statistics of $\gamma$ we remove obvious symmetries (for instance we can assume $x,y\geq 0$) and need only look at those of 
\begin{equation}
\gamma = 16 x^2 y^2 \sin^4 \eta, \;\;\; x,y \;\; {\rm  i.i.d, \;\; with}\;\; \arccos(x)\in_u[0,\pi/2],  \eta\in_u[0,\pi/2]
\end{equation}
In \cite{RealP} the authors gave a complete description of  the totally real subspace in $\IC^2$ of all the parameters for two generator Kleinian groups with one generator elliptic (finite order) of order 2.  Since a group generated by two parabolics admits a $\IZ_2$ extension by adding in the involution $h$ which conjugates $f$ to $g$,  $\langle h,f,g\rangle = \langle f,h \rangle$ and $\langle f,g\rangle \langle f,h\rangle$ with index two both groups are simultaneously discrete or not discrete.  We then deduce the following from \cite{RealP} .
\begin{theorem}\label{FuchsianRiley}  The group $\langle f,g\rangle$ is discrete if and only if $\gamma\geq 16$.
\end{theorem}
Thus we see the p.d.f. for  $\alpha = x y \sin^2 \eta$, with $\arccos(x),\arccos(y)\in_u[0,\pi/2]$, and $\eta\in_u[0,\pi/2]$,  and in the first instance to determine $\Pr\{\alpha \geq 1\}$.  A computational run on $10^7$ matrix pairs gave the  probability $0.3147\ldots$.

The p.d.f. for the random variable $x$ is $F_x=\frac{2}{\pi} \; \frac{1}{1+t^2}$, using the Mellin transform we find the p.d.f. for $xy$ as
\begin{eqnarray*}
F_{xy}(s) & = & \frac{4}{\pi^2} \; \int_{0}^{\infty} \frac{1}{1+t^2} \times \frac{1}{1+s^2/t^2} \times \frac{dt}{t} \\
& = & \frac{4}{\pi^2} \; \frac{\log s}{s^2-1} 
\end{eqnarray*}
a smoothly decreasing function.  Next,  the p.d.f. for $\sin^2\eta$ is 
\[ F_{\sin^2(\eta)} = \frac{1}{\pi} \frac{1}{\sqrt{t(1-t)}},  \hskip20pt 0\leq t \leq 1 \]
Thus we need to compute the Mellin convolution
\begin{eqnarray*}
F_{xy\eta}(s) & = & \frac{4}{\pi^3} \; \int_{0}^{1}   \frac{1}{\sqrt{t(1-t)}} \times \frac{\log s/t}{s^2/t^2-1}  \frac{dt}{t} \\
& = & \frac{4}{\pi^3} \; \int_{0}^{1}    \sqrt{\frac{t}{1-t}} \times \frac{\log s/t}{s^2-t^2} \; dt  
\end{eqnarray*}
Mathematica will do this integral for us,  but the result is rather long and complicated.  We do not reproduce it here,  but have graphed the p.d.f. below.
\begin{center}
\scalebox{0.3}{\includegraphics*[viewport=0 50 800 550]{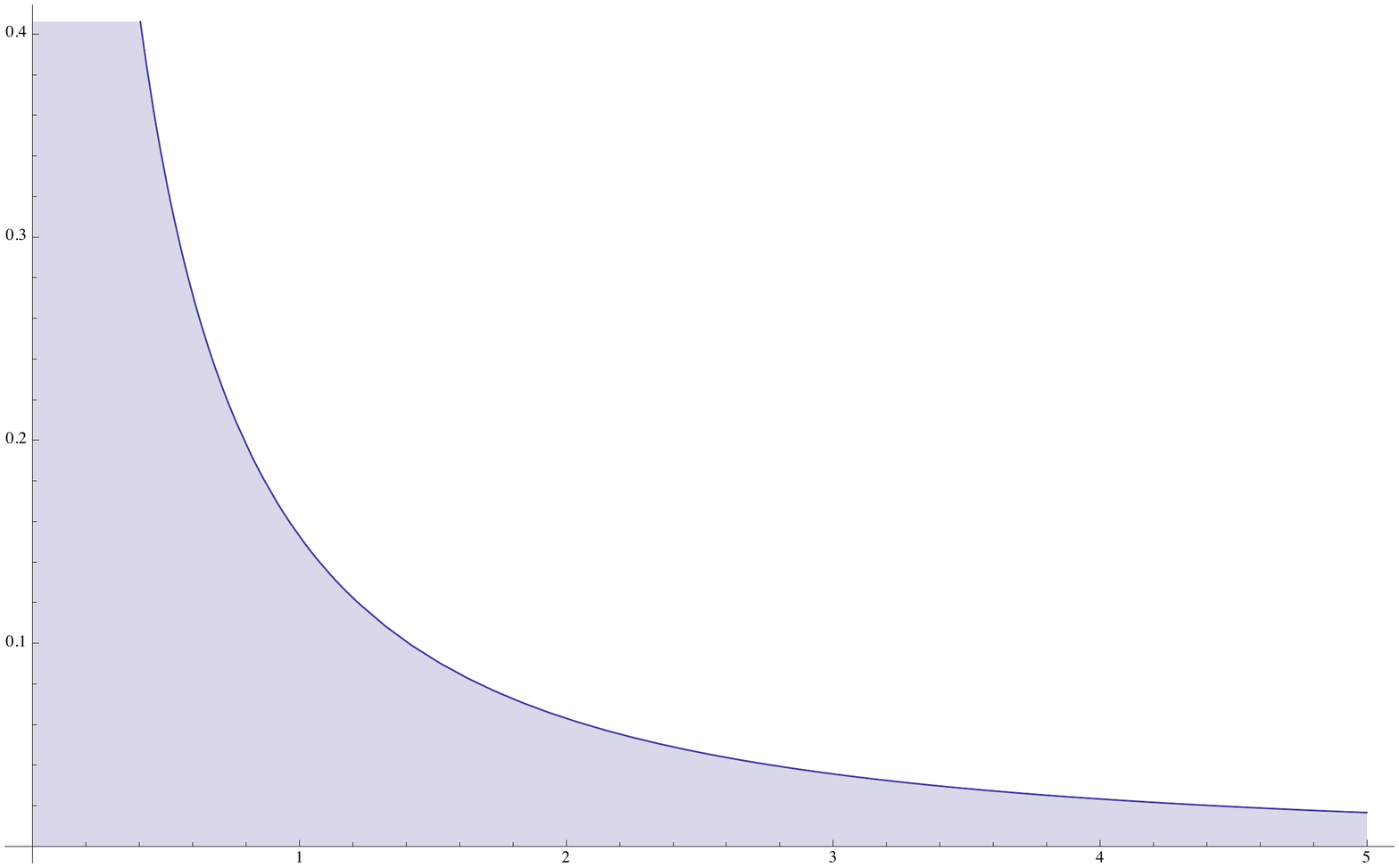}} \\
p.d.f. of  $\sqrt{\gamma/16}$, with $\gamma=\gamma(f,g)$,  $f,g$ parabolic.
\end{center}
What we want is the probability that this number is greater than one.  Thus we compute
\begin{eqnarray*}
\lefteqn{\frac{4}{\pi^3} \; \int_{1}^{\infty} \int_{0}^{1}    \sqrt{\frac{t}{1-t}} \times \frac{\log s/t}{s^2-t^2} \; dt  \; ds}\\ & = & \frac{4}{\pi^3} \; \int_{0}^{1}   \int_{1}^{\infty}   \sqrt{\frac{t}{1-t}} \times \frac{\log s/t}{s^2-t^2} \; ds \; dt \\
 & = & \frac{1}{\pi^3} \; \int_{0}^{1}   \frac{\pi^2-2\log(t) \log \Big(\frac{1+t}{1-t}  \Big) +{\rm PolyLog}\big(2,t^{-2}\big)-4{\rm PolyLog}\big(2,t^{-1}\big) }{ \sqrt{t(1-t)} } \; dt \\
  & = & 1+\frac{1}{\pi^3} \; \int_{0}^{1}   \frac{-2\log(t) \log \Big(\frac{1+t}{1-t}  \Big) +{\rm PolyLog}\big(2,t^{-2}\big)-4{\rm PolyLog}\big(2,t^{-1}\big) }{ \sqrt{t(1-t)} } \; dt 
\end{eqnarray*}
The integral 
\[   \int_{0}^{1}   \frac{{\rm PolyLog}\big(2,t^{-2}\big)-4{\rm PolyLog}\big(2,t^{-1}\big) }{ \sqrt{t(1-t)} } \; dt \]
evaluates in terms of hypergeometic functions to
\[ \pi  \left(\frac{4}{9} \, _4F_3\left(\frac{3}{4},\frac{3}{4},\frac{3}{4},\frac{5}{4};\frac{3}{2},\frac{7}{4},\frac{7}{4};1\right)-8 \, _4F_3\left(\frac{1}{4},\frac{1}{4},\frac{1}{4},\frac{3}{4};\frac{1}{2},\frac{5}{4},\frac{5}{4};1\right)\right) \]
but we are unable to do the integral
\[ \int_{0}^{1}   \frac{\log(t) \arctanh(t) }{ \sqrt{t(1-t)} } \; dt  \approx -0.690591 \] 
but we can evaluate it to high precision.  In this way we find $\Pr\{\alpha \geq 1\} = 0.314833\ldots$ to establish the following theorem.
\begin{theorem}  Let $\Gamma$ be a Fuchsian group generated by two randomly chosen parabolic  elements  of ${\cal F}$.  Then the probability that $\gamma$ is discrete is equal to $0.314833$,
\[ \Pr\{ \Gamma\;\; \mbox{is a discrete Fuchsian group}\} = 0.314833\ldots \]
\end{theorem}
We gave an earlier lower estimate on this probability of $\frac{1}{6}$ using the``ping-pong'' lemma and a discreteness test based on the disjointness of isometric circles.

\begin{theorem}  Let $\Gamma$ be a Fuchsian group generated by two randomly chosen parabolic elements of ${\cal F}$.  Then the probability that the isometric circles of $f$ and $g$ are pairwise disjoint is $\frac{1}{6}$.
\end{theorem}

\section{$\IZ_2$-extensions.}

Every Fuchsian group generated by two parabolics,  say $\Gamma=\langle f,g\rangle$, admits a $\IZ_2$ extension by introducing an elliptic element of order two,  say $\Phi$,  with the property that $\Phi \circ f\circ \Phi^{-1}= g$.  This is discussed in \cite{GehMar} and elsewhere.  The groups $\Gamma$ and $\Gamma_\Phi = \langle f,\Phi\rangle = \langle g , \Phi \rangle = \langle f,g,\Phi \rangle$ are simultaneously discrete or not discrete as $|\Gamma:\Gamma_\Phi|\leq 2$.  Thus we might expect that the groups generated by a randomly chosen parabolic and a randomly chosen elliptic of order two have the same probability of being discrete as that of two randomly chosen parabolics.  We now investigate this.  It is seems clear that conjugacy cannot preserve the uniform distribution of isometric circles,  unless $\Phi$ is an isometry.  What is surprising is that the probability that $\langle f,\Phi \rangle$ is discrete is larger for random parabolic $f$ and elliptic $\Phi$ than the probability that $\langle f,g\rangle$ is discrete.

\medskip

The transformation given at  (\ref{mob})  is elliptic of order two if and only if $\zeta = \pm i$.  Thus we have matrix representations
\begin{equation}\label{2&para} f \leftrightarrow \pm \left( \begin{array}{cc} 1+ix & x e^{i\theta} \\ x e^{-i\theta} & 1-i x \end{array} \right), 
\Phi \leftrightarrow \pm \left( \begin{array}{cc} \frac{i}{\sqrt{1-|w|^2}} & \frac{i\bar w}{\sqrt{1-|w|^2}}\\ \frac{-iw}{\sqrt{1-|w|^2}} & \frac{-i}{\sqrt{1-|w|^2}}  \end{array} \right),  \end{equation}
with $x,w$ chosen as described above; $x=\cot(\eta)$,  $|w|=\cos(\alpha)$ and the angles $\eta,\alpha$ uniformly distributed in $[0,\pi/2]$.

Then if we define 
$g = \Phi f  \Phi^{-1}$
we compute  that
\begin{eqnarray*}
g & = & \left(
\begin{array}{cc}
 1+\frac{i x\left(|w|^2-2 \sin (\alpha +\theta ) |w|+1\right)}{1- |w|^2} & \frac{e^{-i (2 \alpha +\theta )} \left(i e^{i (\alpha +\theta )}+ |w|\right)^2 x}{1- |w|^2} \\
 \frac{e^{i (2 \alpha +\theta )} \left(-i e^{-i (\alpha +\theta )}+ |w|\right)^2 x}{1- |w|^2} & 1-\frac{i x \left( |w|^2-2 \sin (\alpha +\theta ) |w|+1\right)}{1- |w|^2} \\
\end{array}
\right)
\end{eqnarray*}
The fixed point of $g$ is  \[ \zeta=e^{-i \alpha }  \; \frac{  i t-e^{i (\alpha +\theta )} }{t e^{i (\alpha +\theta )}-i} \]
and it is straightforward to see this is uniformly distributed on the circle.   We are therefore led to consider the distribution of the arclength of isometric circles,
\begin{eqnarray*} 
\cot^{-1}\Big(\frac{x \left(1+|w|^2-2 |w| \sin (\alpha +\theta )\right)}{1-|w|^2}\Big) =\cot^{-1}\Big(   \frac{\cot(\eta) \left(1+\cos^2(\alpha)-2\cos(\alpha) \sin (\alpha +\theta )\right)}{\sin^2(\alpha)}\Big).
\end{eqnarray*}
This distribution is not uniform.  A p.d.f for all of the individual terms in this expression can be found,  but the expression for the p.d.f. of this random variable is very complicated. We simply generate the histogram shown below.

\begin{center}
\scalebox{0.7}{\includegraphics*[viewport=60 470 600 800]{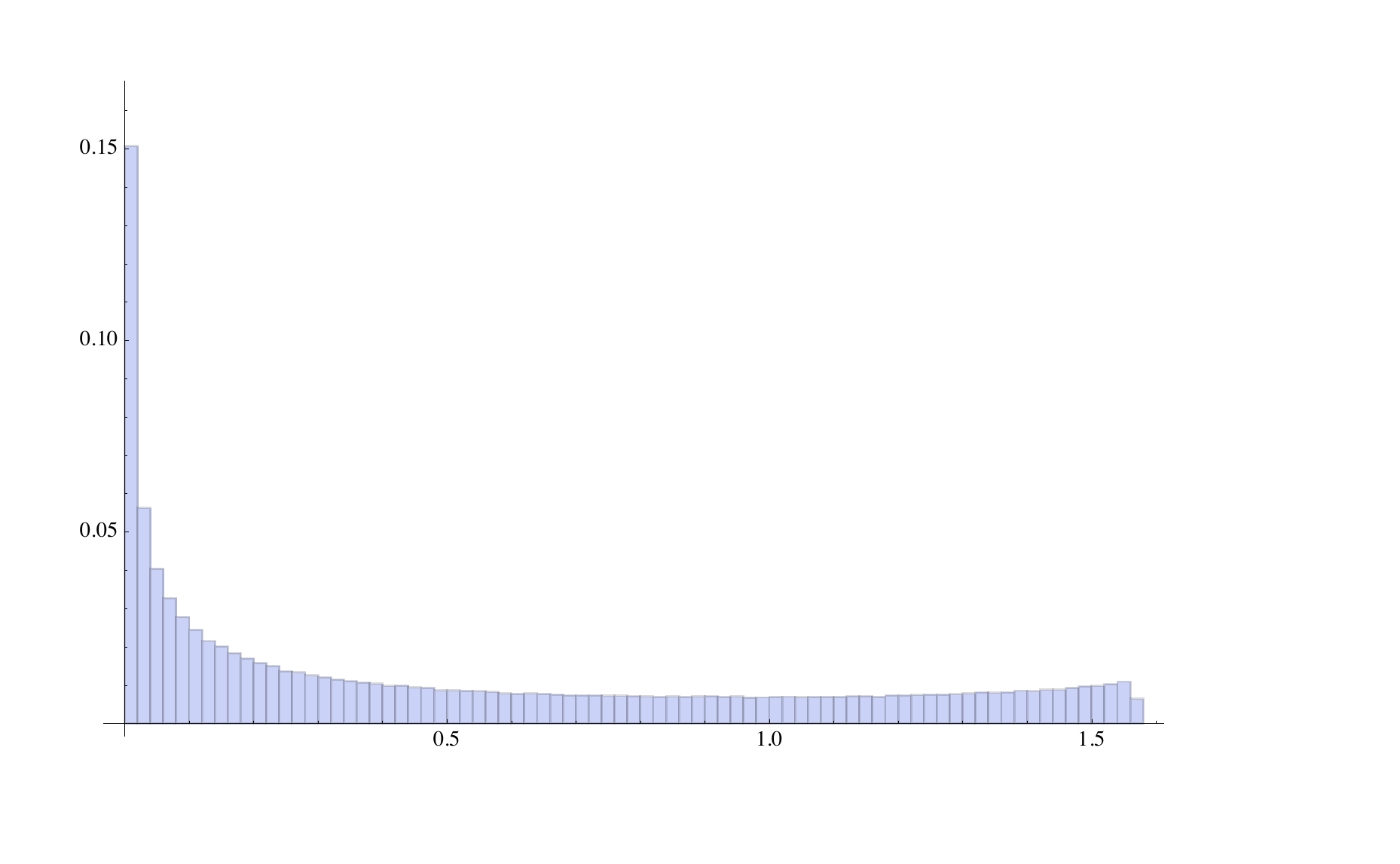}} \\
$\cot^{-1}\Big(   \frac{\cot(\eta) \left(1+\cos^2(\alpha)-2\cos(\alpha) \sin (\alpha +\theta )\right)}{\sin^2(\alpha)}\Big)$,  $\eta,\alpha\in_u[0,\pi/2]$ and $\theta\in_u[0,2\pi]$
\end{center}
This histogram shows that the distribution of the arclength of isometric circles slightly favours smaller arclengths.  Since disjointness of these circles is predicts discreteness we expect that the group $\langle f,g \rangle$ is more likely discrete than that of a group generated by two random parabolics.

The commutator parameter is 
\[ \gamma(f,\Phi) = \frac{4 x^2 (1-|w| \sin (\alpha +\theta ))^2}{1-|w|^2}  = \sqrt{\gamma(f,g) }\]
The group $\Gamma_\Phi$ is discrete and free if and only if $\gamma(f,\Phi) \geq 4$,  yielding the inequality
$\cot(\eta) (1-\cos(\alpha) \sin (\alpha +\theta )) \geq \sin(\alpha)$.
 Thus we want to know the probability that 
 \[ \sec (\alpha )-\tan (\alpha ) \tan (\eta )>\cos (\theta ).\] 
 An experiment on $10^7$ triples of random variables shows this to happen with probability about $0.595$.   We were unable to identify the p.d.f. to get an explicit integral,  but it is not too difficult to achieve the following computationally supported result.
 \begin{theorem}  Let $\Gamma=\langle f,\Phi \rangle$ be a Fuchsian group generated by a random parabolic $f$ and a random element of order two $\Phi$.  Then the probability that $\Gamma$ is discrete is $a_{\infty,2}$,  where
 \[ a_{\infty,2} =\frac{4}{\pi^3}\; \Big|\big\{(\theta,\alpha,\eta)\in [0,2\pi]\times [0,\pi/2]\times [0,\pi/2] : \sec (\alpha)-\tan (\alpha ) \tan (\eta )\geq \cos (\theta ) \big\}\Big| \]
 and 
 \begin{equation}
 0.59 \leq a_{\infty,2} \leq 0.6.
 \end{equation}
 \end{theorem}
 \noindent{\bf Proof.}  We rewrite the inequality as 
 \[ F(\theta,\alpha,\eta) =\cos(\eta) \big(1-\cos(\alpha) \sin (\theta )\big) -\sin(\eta) \sin(\alpha)  \geq 0.\]
 We then compute that
\begin{eqnarray*} |\nabla F| & = &  |F_\theta|+|F_\eta|+|F_\alpha| \\& = &  |\cos (\alpha ) \cos (\eta ) \cos (\theta )|+|\sin (\alpha ) \cos (\eta ) \sin (\theta )-\cos (\alpha ) \sin (\eta )| \\ && +|\sin (\eta ) (\cos (\alpha ) \sin (\theta )-1)-\sin (\alpha ) \cos (\eta )| \\
& \leq &  2+2 |\cos (\eta+\alpha )| \leq 4 
\end{eqnarray*}
We partition the region $[0,2\pi]\times [0,\pi/2]\times [0,\pi/2]$ initially into approximately 2 million boxes of side length $\frac{1}{50}$.  We use a first order Taylor estimate.  We evaluate $F$ at each vertex (relying on Mathematica to compute the value accurately to 4 decimal places).  If at each vertex $v_i$ we have $F(v_i)-4\sqrt{3}/100\geq 0.001$ we admit the box.  This is the coarse first pass.  When a box is not admited we subdivide each edge,  obtain a better gradient bound $M=|\nabla F|_Q$ on the box,  and apply the test  $F(v_i)-M\sqrt{3}/200\geq 0.001$ and admit the subdivided box or not.  We only record a count of the boxes admitted (and if it is a subdivided box) and subdivide at most twice.  From this we get a count of the number of boxes admitted  and their volume which we then sum.  This number is reported as the lower bound.  An upper bound is computed identically but the test is less than or equals and we use $-0.001$.  A number of easy simplifications can be made.  For instance if  $\sin(\theta)\leq 0$,  then the test is implied by
\[ \cos(\eta) - \sin(\eta) \sin(\alpha)  \geq 0 \]
removing $0\leq \eta \leq \pi/4$ and all $\alpha$,  a quarter of the search space. \hfill $\Box$

\medskip

The statement of the proof infers the existence of the numbers $a_{\infty,m}$ for the probability a group generated by a random parabolic and element of order $m$.  The story here is a little more complicated and will be discussed elsewhere.  The discreteness bounds on the commutator are also different. 

\section{Kleinian Groups.}

We take our cue from the Fuchsian case regarding the definition of a random parabolic element where we discovered that the isometric disks met the unit circle in uniformly distributed pairs of arcs.  We want a probability distribution on the parabolic subgroup
\[ {\cal P} = \{ A \in SL(2,\IC): \tr^2 A - 4 = 0 \} \]
We seek a distribution invariant under the conjugation action by spherical isometries of the Riemann sphere $\oC$.  These spherical isometries are represented by the matrices (ignoring the  $\pm $  term)
\begin{equation}\label{si}
I = \left( \begin{array}{cc} a & b \\ c & d \end{array} \right), \;\; |a|^2+|b|^2+|c|^2+|d|^2 = 2, \;\; ad - bc = 1.\;\; \leftrightarrow \frac{az+b}{cz+d}
\end{equation}
(see \cite[]{Beardon}).  What we mean here is that given a set $U\subset {\cal P}$ and a randomly chosen $f$ (actually we choose $\langle f\rangle$) we have for any spherical isometry $\phi$
\[ {\rm Pr} \{ f\in U \} = {\rm Pr} \{ f\in  \phi \circ U \circ \phi^{-1} \} \]
where $\phi \circ U \circ \phi^{-1}= \{\phi \circ f \circ \phi^{-1}: f\in U\}$.  Notice that $\phi {\cal P} \phi^{-1} = {\cal P}$.
Those parabolic M\"obius transformation fixing $\infty\in \oC$ are easy to describe.  They have the form $z\mapsto z+\lambda e^{i\theta}$ with $\lambda >0$ and $\theta\in [0,2\pi)$.  The isometric circles of Parabolic transformations whose fixed point is not infinity are tangent at the parabolic fixed point. We could now look at the intersections of the isometric spheres of the Poincar\'e extensions (to $\hat{\IR^3}$) with $\oC$ and create a distribution on the entries so that the solid angles are uniformly distributed.  This is quite complicated, and we prefer the equivalent intrinsic approach following.  

\medskip

We choose the parabolic fixed point $z_0$ uniformly in the spherical measure of $\oC \approx \IS^2$.  This measure comes from the spherical {\em metric} on $\oC$ defined by
\[ ds^2 = \frac{2}{1+|z|^2} \]

We next choose a solid angle $\eta\in [0,\pi/2]$ uniformly $\eta$ is the spherical radius of a spherical disk in $\oC$.  We now want to identify 
a pair of spherical disks $D_1$ and $D_2$ of radius $\eta$ mutually tangent to $z_0$.  The parabolic M\"obius transformation $f$ we seek will have $f(D_1)=\oC\setminus D_2$ and $f(z_0)=z_0$.  Of course the order of these disks is important here for the M\"obius transformation that has $h(D_2)=\oC\setminus D_1$ and $h(z_0)=z_0$ is easily found to be $h=f^{-1}$.  However this means $\langle f \rangle = \langle h \rangle$ and this is all we care about.  There are a ``circles worth'' of choices for this pair of disks.  Thus we choose randomly and uniformly from the family of great circles through $z_0$ (which we identify with the uniform distribution on $\IS$).  Call this circle $C$,   and then choose the pair of spherical disks $D_1$ and $D_2$ of radius $\eta$ mutually tangent to $z_0$ so that the (spherical) centers of these disks (and $z_0$) lie on $C$.  It is immediate that this selection process is invariant under the conjugation action of the spherical isometry group as all the measures are.  Thus we can now go about the task of identifying the matrix entries of a randomly selected cyclic parabolic group (we abuse notation and say ``random parabolic'').

\medskip

Given a pair of spherical disks of solid angle $\eta$ tangent at $\infty \in \oC$ and whose centers are real,  the associated parabolic transformation of $\oC$ has the form $z\mapsto z+\lambda$,  $\lambda >0$.  The number $\lambda$ can be found as the Euclidean distance in $\IC$ between intervals $\oR\cap D_1$ and $\oR \cap D_2$.  The spherical length of $\oR$ is $2\pi$,  and so the spherical distance between these intervals is $2\pi -4\eta$.  The midpoint of this gap is $0$ and so we have
\begin{equation}
\pi - 2\eta = \int_{0}^{\lambda/2} \frac{|dz|}{1+|z|^2}=2\tan^{-1}\Big(\frac{\lambda}{2}\Big)
\end{equation}
and,  as with the Fuchsian case,  we see $\lambda = \cot(\eta)$.  The selection of a great circle through which the centers of these disks lie is effected by replacing $\lambda$ by $e^{i\theta}\lambda$,  with $\theta\in_u[0,2\pi]$. This has given us the matrix presentation
\begin{equation}\label{trans}
P = \left( \begin{array}{cc} 1& e^{i\theta}\lambda \\ 0 &1\end{array}\right) 
\end{equation}
We now need to identify any (!) spherical isometry which sends $z_0$ to $\infty$.  

The spherical isometries which stabilise the extended real axis have the matrix form
\[ \Phi = \left( \begin{array}{cc} \cos(\theta) & \sin(\theta) \\-\sin(\theta) & \cos(\theta) \end{array}\right) \]
So $\Phi(|z_0|)=\infty$ requires $\sin(\theta)|z_0|=\cos(\theta)$.  That is $ \theta = \cot^{-1}(|z_0|) $ with the matrix representation
\begin{equation} \Phi = \left( \begin{array}{cc} \frac{|z_0|}{\sqrt{1+|z_0|^2}} & \frac{1}{\sqrt{1+|z_0|^2}} \\\frac{-1}{\sqrt{1+|z_0|^2}}  & \frac{|z_0|}{\sqrt{1+|z_0|^2}}\end{array}\right) 
\end{equation}
Then $z_0\mapsto |z_0|$ is effected by the spherical isometry with matrix representation
\[ \left( \begin{array}{cc}  e^{-i \arg(z_0)/2} & 0 \\0 &e^{i \arg(z_0)/2}  \end{array}\right) \]
Therefore a spherical isometry we could choose has the form (with $\eta=\arg(z_0)/2$)
\[ \Psi = \left( \begin{array}{cc} \frac{|z_0|e^{-i \eta} }{\sqrt{1+|z_0|^2}} & \frac{e^{i \eta}}{\sqrt{1+|z_0|^2}} \\\frac{-e^{-i \eta}}{\sqrt{1+|z_0|^2}}  & \frac{|z_0|e^{i\eta} }{\sqrt{1+|z_0|^2}}\end{array}\right) \]
Thus the parabolic element associated with a random pair of disks tangent at $z_0$ has the form
\begin{eqnarray*}
\Psi^{-1} \circ P \circ \Psi & = & \ \left(
\begin{array}{cc}
 1-\frac{e^{i \theta } t \lambda }{t^2+1} & \frac{e^{i (\eta +\theta )} t^2 \lambda }{t^2+1} \\
 -\frac{e^{-i (\eta -\theta )} \lambda }{t^2+1} & \frac{e^{i \theta } t \lambda }{t^2+1}+1 \\
\end{array}
\right), \hskip10pt t=|z_0| \\
& = & \II + \frac{e^{i \theta }  \lambda }{t^2+1}  \left(
\begin{array}{cc}
-  t    &  e^{i  \eta } t^2  \\
 - e^{-i \eta }   &   t   \\
\end{array}
\right), \hskip10pt t=|z_0| 
\end{eqnarray*}
where $\II$ is the identity matrix. If we set for $i=1,2$
\[ F_i =  \II + \frac{e^{i \theta_i}  \lambda_i}{t_i^2+1}  \left(
\begin{array}{cc}
-  t_i    &  e^{i  \eta_i } t_i^2  \\
 - e^{-i \eta_i }   &   t_i   \\
\end{array}
\right), \hskip10pt t_i=|z_i|  \]
and put $G=[F_1,F_2]$ Then we may compute that 
\begin{eqnarray*}\gamma(F_1,F_2) & = &  \tr(G)-2 = \frac{\lambda_1^2 \lambda_2^2 e^{2 i (\theta_1+\theta_2)} \left(t_1 e^{\frac{1}{2} i (\eta_1-\eta_2)}-t_2 e^{-\frac{1}{2} i (\eta_1-\eta_2)}\right)^4}{\left(t_1^2+1\right)^2 \left(t_2^2+1\right)^2}\\
& = & \frac{\lambda_1^2 \lambda_2^2 e^{2 i (\theta_1+\theta_2)} (\cos (\eta )  (t_1-t_2)-i \sin (\eta ) (t_1+t_2))^4}{\left(t_1^2+1\right)^2 \left(t_2^2+1\right)^2}, \;\; \eta=\frac{\eta_1-\eta_2}{2}\\
\end{eqnarray*}
We are interested in the distribution of the number $\gamma=\gamma(F_1,F_2)$.  Since  $\theta_1,\theta_2\in_u [0,2\pi]$ we may replace the sum $\theta_1+\theta_2=\theta\in_u[0,2\pi]$.  It is also clear from this form that the argument of $\gamma$ is uniformly distributed.  We next see that
\[ |\gamma| =  \frac{\lambda_1^2 \lambda_2^2  \big(t_1^2-2 t_1t_2 \cos (2 \eta )+t_2^2\big)^2}{\left(t_1^2+1\right)^2 \left(t_2^2+1\right)^2} \]
We recall that
\begin{eqnarray*}
t_i = \tan\big(\alpha_i \big) && \alpha_i\in_u [0,\pi/2] \\
\lambda_i = \cot \big(\beta_i\big) && \beta_i\in_u [0,\pi/2]
\end{eqnarray*}
so $t_i$ and $\lambda_i$ are identically distributed, $i=1,2$ and hence 
\[ |\gamma| =  \cot^2(\beta_1) \cot^2(\beta_2) \big[\sin^2(\eta ) \sin^2 (\alpha_1)+\cos^2 (\eta ) \sin^2 (\alpha_2)\big]^2  \]
where we have chosen $\eta\in_u[0,\pi]$ and used the periodicity and symmetry to replace $\alpha_1+\alpha_2$ by $\alpha_1$ and $\alpha_1-\alpha_2$ by $\alpha_2$ without affecting the distribution.

By way of comparison with the Fuchsian case,  we present a similar histogram for the distribution of $|\gamma(f,g)|$ derived from $10^6$ random pairs of parabolic elements generated by Mathematica code.  It is clearly more heavily weighted toward the value $0$.

\begin{center}
\scalebox{0.75}{\includegraphics*[viewport=60 470 600 800]{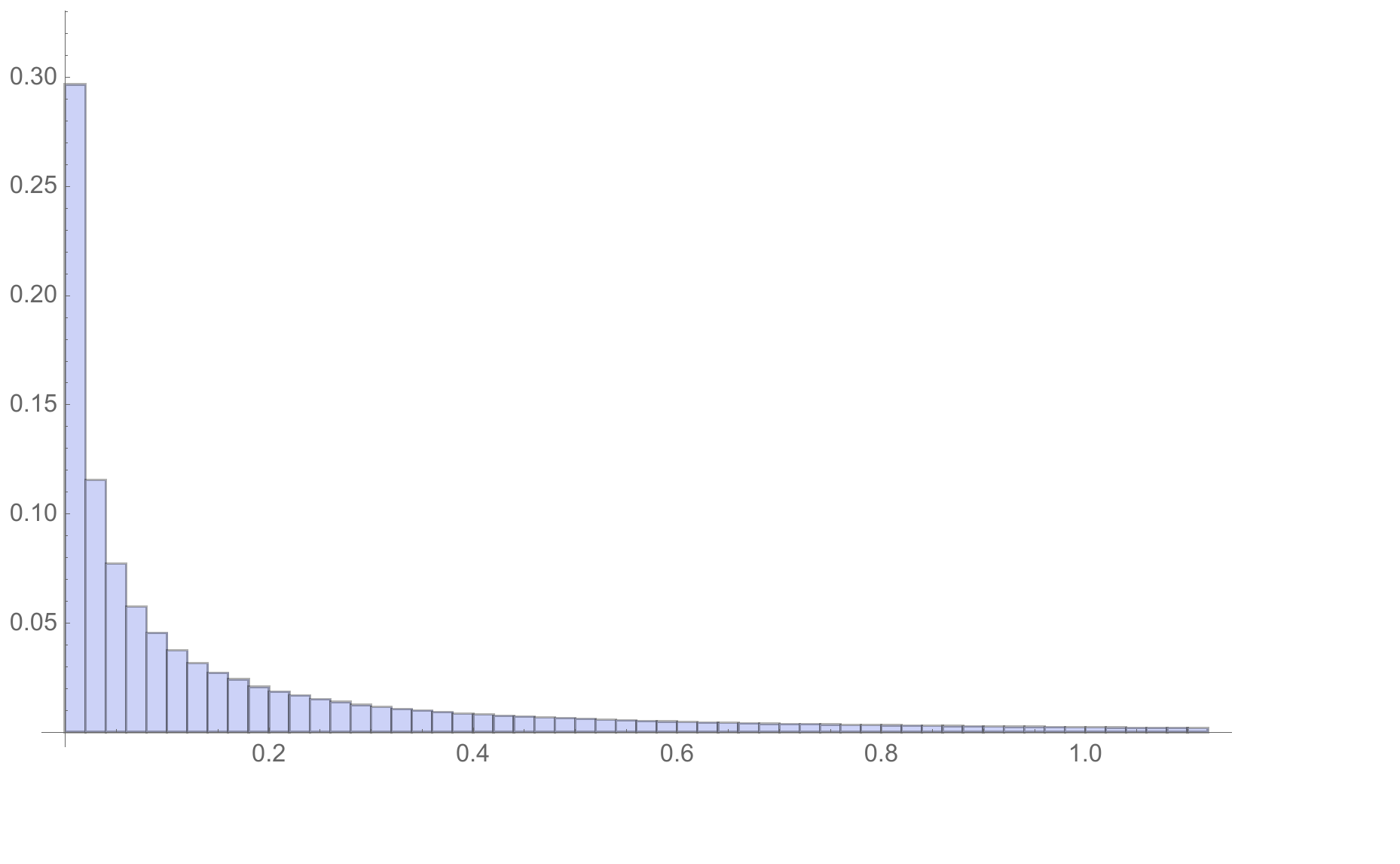}} \\
p.d.f. of  $\sqrt{|\gamma|/16}$, with $\gamma=\gamma(f,g)$,  $f,g$ parabolic.
\end{center}

\subsection{Some common normalisations.}
Discrete groups generated by two parabolics have been studied before and  a common normalisation is that the fixed points are $0$ and $\infty$.  In general fixed points will not be antipodal,  the expected angle between them is in fact $\pi/2$.  However with this conditioning we have two random parabolics have the form 
\begin{equation} \label{PQ}
P = \left( \begin{array}{cc} 1& e^{i\theta_1}\lambda_1 \\ 0 &1\end{array}\right), \;\;\;\;  Q = \left( \begin{array}{cc} 1& 0 \\ e^{i\theta_2}\lambda_2  &1\end{array}\right).
\end{equation}
Then,  with the same distributions as above,
 \begin{eqnarray*}
\gamma(P,Q)  &=& e^{2i\theta}  \cot^2(\beta_1) \cot^2(\beta_2)  
\end{eqnarray*}
This is the same as the distribution for $\frac{1}{z^2} \cdot \frac{1}{w^2} $ for random $z,w\in \oC$ with the spherical metric,  and since inversion is a spherical isometry this has the same $p.d.f.$ as $z^2 \cdot w^2$.

\medskip
Another common nomalisation is to assume that in (\ref{PQ}) we have $Q_{2,1}=1$.  We then define
\begin{equation}\label{gammau}
 \Gamma_u=  \langle\; \left( \begin{array}{cc} 1&u\\ 0 &1\end{array}\right),  \left( \begin{array}{cc} 1& 0 \\ 1 &1\end{array}\right)\rangle 
 \end{equation}
and construct the set
\begin{equation}
{\cal R} = \{u\in \IC: \mbox{the group $\Gamma_u$ is discrete and freely generated} \},
\end{equation}
is called the Riley slice. It is closed,  connected and the boundary is a topological circle \cite{}. Since the commutator parameter $\gamma_G=\gamma(A,B)$ where $G=\langle A,B\rangle$ is a conjugacy invariant,  and since two parabolics $A$ and $B$ can be simultaneously conjugated to the form $P$ and $Q$ we have the following theorem which is the analogue in the Kleinian case of Theorem \ref{FuchsianRiley}.
\begin{center}
\scalebox{0.3}{\includegraphics*[viewport=60 -50 700 550]{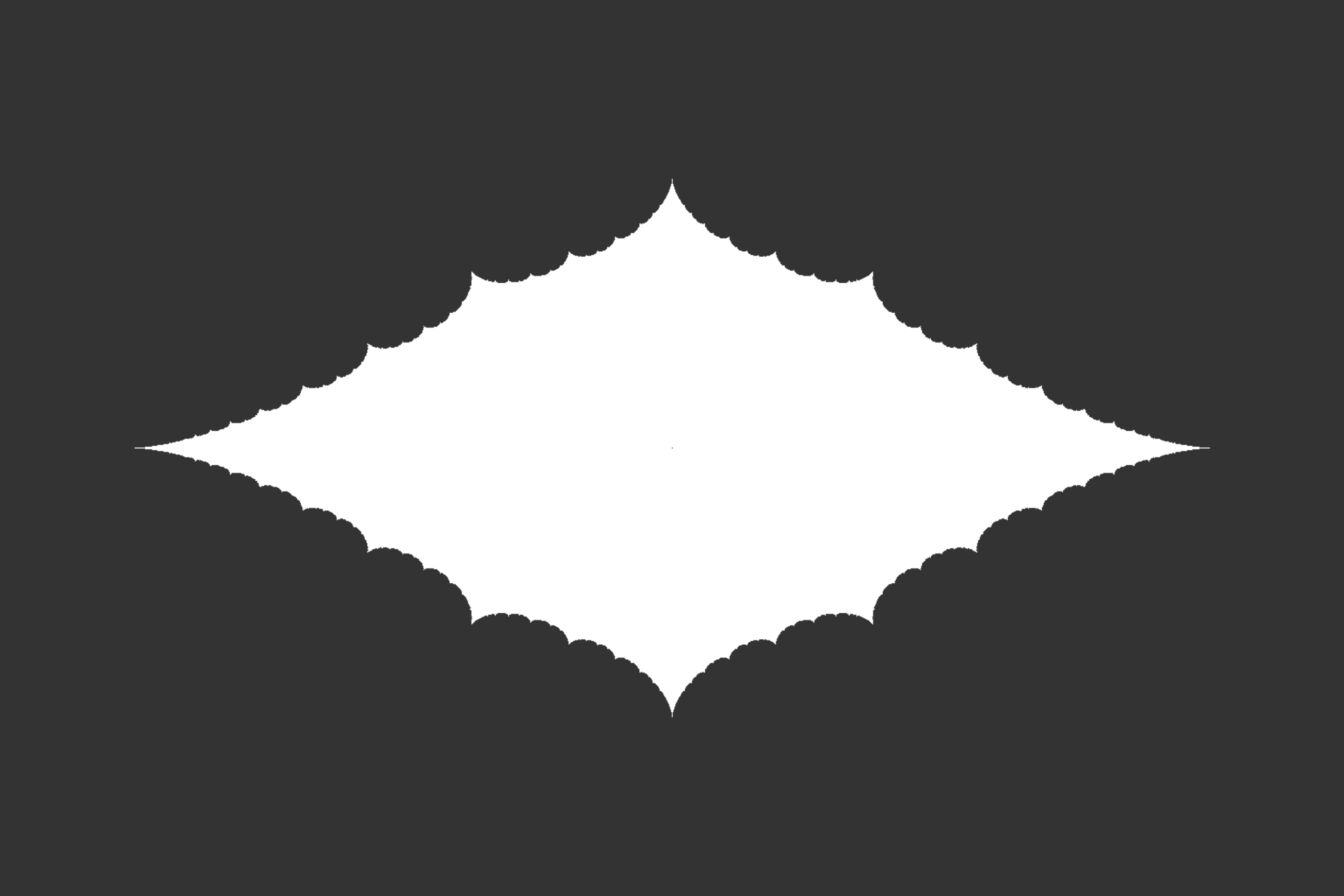}} \\
The Riley Slice (in black)
\end{center}
\begin{theorem}\label{thm5.8}
Let $G=\langle A,B\rangle$,  where $A$ and $B$ are parabolic elements.  Then $G$ is discrete and free if and only if $\sqrt{\gamma} \in R$.
\end{theorem}
Since $R$ is symmetric,  the branch of square root is immaterial.  However the Riley slice is a rather complicated set and there is little chance of providing any analytical results.  We also note that the space of discrete groups which are not free has measure $0$,  in fact consists of a discrete set accumulating on the boundary of the Riley slice so from a probabilistic point of view we can ignore them.

\subsection{An experimental approach.}

In \cite{Bow} B. Bowditch studies the behaviour of trees of Markoff triples defined over the complex numbers and relates this to the space of type-preserving representations of the punctured torus group into $SL(2,\IC)$. In particular, he considers Markoff triples correspond to quasifuchsian representations. Of relevance to us,  he derives a variation of McShane's identity for quasifuchsian groups and he relates,  in the case of non-discrete representations, the asymptotic behaviour of Markoff triples and the realisability of laminations in hyperbolic 3-space.  For us,  this gives a conjectural description of the exterior of the Riley slice to study computationally.  
In particular the representation at (\ref{gammau}) corresponds to the (type-preserving) Markoff map $\phi $ defined by
 $\phi(0/1) = x$,  $ \phi(1/0) = 0$,  and  $\phi(1/1) = x$
where $u = -x^2$.
(For example, \cite[Lemma 5.3.2]{ASWY}).
The (type-preserving) Markoff map $\phi$ is valued in the complement of Bowditch's space
if there exists $p/q$ such that $|\phi(p/q)| < 1$.
(an unpublished extension of \cite[Theorem 1.1]{NT}).   We can use $\phi(p/q)$ for any $p/q\in \IQ$
to get a description of (the interior of) the complement of Bowditch's space.  The algorithmic study of this conjecture and related questions is also addressed in \cite{SPY}.   Using these ideas we can estimate the spherical area of the exterior of the Riley slice to be $0.779\times 4\pi$.   

\medskip

We can use this as follows.  The group $\Gamma_{\sqrt{\gamma}}$ is discrete if and only if $\sqrt{\gamma}\in {\cal R}$.  From our earlier calculation we see that $\sqrt{\gamma}$ is uniformly distributed in the spherical metric when the parabolic generators are normalised as per the definition of $\Gamma_u$.  We  hence observe the following alignment of our definitions.

\begin{lemma}  If $u$ is uniformly selected in the spherical metric,  then $\Gamma_u$ is a group generated by a random parabolic $f$ with fixed point $0$ together with the parabolic $g(z)=z+1$.  \end{lemma}

We therefore estimate that the probability a random $\Gamma_u$ is not discrete to be around $0.779$.  However,  as the commutator is a conjugacy invariant we can use Theorem \ref{thm5.8} to estimate the probability of non-discreteness in the general case. Using the description of the complement of the Riley slice as per Bowditch's condition as described above we can run a trial to determine a rough estimate of the probability that 
\[ e^{i\theta}\cot(\beta_1) \cot(\beta_2) \big[\sin^2(\eta ) \sin^2 (\alpha_1)+\cos^2 (\eta ) \sin^2 (\alpha_2)]\in {\cal R},\]
where $\eta, \beta_i,\alpha_i\in_u[0,\pi/2], \;\;\theta\in_u[0,2\pi] $,  and hence the probability that a group generated by two random parabolics is discrete.  This is described in the following figure and from this we estimate that the the probability of a group generated by two random parabolics is not discrete is approximately $0.768$.  
 
\begin{center}
\scalebox{0.5}{\includegraphics*[viewport=-30 -50 550 550]{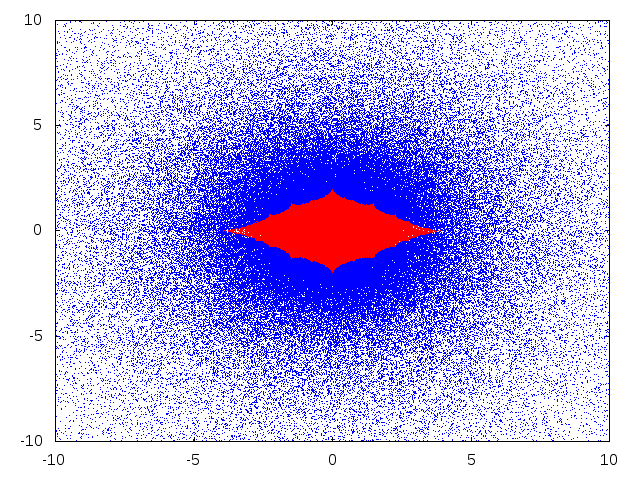}} \\
Bowditch condition for non-discrete groups in red.  
\end{center}
The figure also shows $10^6$ commutators of randomly selected parabolic elements in blue.  Approximately  $76.8\%$ of these commutators lie in the red region and thus cannot be those of a discrete group.

\medskip

\noindent {\bf Conjecture:}  The probability that a group generated by two random parabolic elements is discrete is less than $0.232$.

\subsection{Estimates on discreteness.}

We can support these calculations with a provable result.  However,  it is quite clear we are a long way form what might be sharp.  First,  the Shimitzu-Leutbecher inequality \cite{Beardon},  basically a version of J\o rgensen's inequality for parabolic generators,   identifies the unit disk as being outside the Riley slice (so in the red region above).  This has the implication that any group $\langle f,g \rangle$ with $f$ parabolic is not discrete if $|\gamma(f,g)|<1$ or unless the group is a Euclidean triangle group (and $\gamma(f,g)=0$).  Since $\gamma(f,g)=0$  happens with probability zero  a random group generated by two parabolic elements is not discrete if 
\[ 0<   \cot(\beta_1) \cot(\beta_2) \big[\sin^2(\eta ) \sin^2 (\alpha_1)+\cos^2 (\eta ) \sin^2 (\alpha_2)\big] < 1  \] 
with the angular distributions as above. We observe that 
\[  \sin^2(\eta ) \sin^2 (\alpha_1)+\cos^2 (\eta ) \sin^2 (\alpha_2) \leq 1 \]
Thus $|\gamma(f,g)| <1$,  and nondiscreteness,  is implied by
\[ |\gamma(f,g)| \leq \cot(\beta_1) \cot(\beta_2). \]
For $\beta_i\in_u[0,\pi/2]$,  $\cot(\beta_1) \cot(\beta_2) <1$ implies $\cos(\beta_1) \cos(\beta_2) - \sin(\beta_1)\sin(\beta_2)< 0$.  That is $\cos(\beta_1+\beta_2)<0$, which is equivalent to $\beta_1+\beta_2>\pi/2$ and this probability is clearly equal to $\frac{1}{2}$.

Next we consider the question of when the   isometric circles are disjoint. It is an elementary consequence of the Klein combination theorem,  sometimes called the ``ping-pong'' lemma,  that the group $G=\langle f,g\rangle$ generated by two random parabolics will be discrete if the isometric circles are disjoint - recall here we are considering the isomeric circles being paired in the spherical metric,  these will not in general be the same as in the euclidean metric of $\oC$.   Let $z$ and $w$ be the fixed points of $f$ and $g$ respectively,  with isometric circles of solid angle $\eta_f$ and $\eta_g$.  If the spherical disks of solid angle $2\eta_f,2\eta_g \in_u [0,\pi]$ are disjoint,  then certainly the isometric circles are disjoint. Since $z$ and $w$ lie on a unique great circle with probability $1$,  and the angle,  say $\theta$,  between $z$ and $w$ on this great circle is uniformly distributed in $[0,\pi]$ we arrive at the question of computing the probability that $\theta_1+\theta_2<\theta_3$ for $\theta_i\in [0,\pi]$.  We leave it to the reader to determine that this probability is $\frac{1}{6}$ - we did this for the Fuchsian case earlier.  Taken together these results give a proof for the next theorem.

\begin{theorem}  Let $G=\langle f,g \rangle$ be generated by two random parabolics.  Then the probability that $G$ is discrete is greater than $\frac{1}{6}$ and less than $\frac{1}{2}$.
\end{theorem}

\end{document}